\documentclass[12pt,letterpaper]{amsart}
\usepackage{amssymb}
\usepackage{amsfonts}
\usepackage{amsmath}
\usepackage[matrix,arrow,graph]{xy}

\theoremstyle{plain}
\newtheorem{thm}{Theorem}

\newtheorem{prop}[thm]{Proposition}

\theoremstyle{definition}

\newtheorem*{defn}{Definition}

\newtheorem*{ex}{Example}

\newtheorem*{rmk}{Remark}

\newcommand{\ra}{\mathop{\rightarrow}\limits}

\newcommand{\ol}{\overline}

\def\ker{\mathop{\rm ker}\nolimits}

\def\id{{\rm id}}

\newcommand{\C}{{\mathbb C}}
\newcommand{\Z}{{\mathbb Z}}
\newcommand{\R}{{\mathbb R}}

\newcommand{\Q}{{\mathbb Q}}

\newcommand{\HH}{{\mathbb H}}
\newcommand{\cal}{\mathcal}

\newcommand{\cF}{{\cal F}}

\newcommand{\fl}{\mbox{\rm fl}}
\newcommand{\flag}{\cal F}
\newcommand{\idperm}{1}

\newcommand{\mfS}{{\mathfrak S}}
\begin{document}

\bibliographystyle{siam}

\title[Lower bounds for KL polynomials from patterns]{Lower bounds for 
Kazhdan-Lusztig polynomials from patterns}
\date{\today}
\author{Sara C. Billey}
\address{
Dept.\ of Mathematics, 2-334\\ 
Massachusetts Institute of Technology\\ Cambridge, MA 02139}
\email{billey@math.mit.edu} 

\author{Tom Braden}
\address{
Dept.\ of Mathematics and Statistics\\
University of Massachusetts\\
Amherst, MA 01003}
\email{braden@math.umass.edu}

\thanks{The first author was supported by National Science
Foundation grant DMS-9983797; the second author was supported in
part by grant DMS-0201823 }

\begin{abstract}
Kazhdan-Lusztig polynomials $P_{x,w}(q)$ play an important role in the
study of Schubert varieties as well as the representation theory of
semisimple Lie algebras.
We give a lower bound for the values $P_{x,w}(1)$ in terms of
``patterns''.  A pattern for an element of a Weyl group is its image 
under a combinatorially defined map to a subgroup generated by 
reflections. This  generalizes the classical definition of patterns 
in symmetric groups.  This map corresponds geometrically to 
restriction to the fixed point set of an action of a one-dimensional torus on
the flag variety of a semisimple group $G$.   
Our lower bound comes from applying a 
decomposition theorem for ``hyperbolic localization'' \cite{Br}
to this torus action.
This gives a geometric explanation for the appearance of pattern 
avoidance in the study of singularities of 
Schubert varieties. 
\end{abstract}


\maketitle

\section{Introduction}

Many recent results on the singularities of Schubert varieties $X_w$
in the variety $\flag_n$ of flags in $\C^n$ are expressed by the
existence of certain patterns in the indexing permutation $w \in \mfS_n$.
For example,
Lakshmibai and Sandhya \cite{LS} proved that $X_w$ is singular if and
only if $w$ contains either of the patterns $4231$ or $3412$
(see also \cite{R}, \cite{W}).
A permutation $w\in \mfS_n$ is said to contain the pattern 
$\tilde{w} \in \mfS_k$ for $k < n$ if the permutation matrix of $w$ 
has the permutation matrix of $\tilde{w}$ as a submatrix.

This implies that
if $\tilde{w} \in \mfS_k$ is any pattern for $w$ and
$X_{\tilde{w}} \subset \flag_k$
is singular, then $X_{w}$ is singular as well.
In this paper, we give a general geometric explanation of
this phenomenon which works for the 
flag variety $\flag$ and Weyl group $W$ of any semisimple algebraic 
group $G$.

Our result concerns the Kazhdan-Lusztig polynomials $P_{x,w}(q) \in
\Z_{\ge 0}[q]$, $x,w \in W$.  Although defined purely combinatorially,
they carry important information about representation theory of
Hecke algebras and Lie algebras 
(see \cite{KL1,BeBe,BryK,BGS} among many others), as well as geometric 
information about the singularities of Schubert varieties $X_w$ in 
$\cF$.

More precisely, $P_{x,w}(q)$ is the Poincar\'e polynomial
(in $q^{1/2}$) of the local intersection cohomology of $X_w$ at a 
generic point of $X_x$, and 
 $P_{\idperm,w}(1)=1$ if and only if $X_{w}$ is rationally smooth \cite{KL2}.
If $G$ is of type A,D, or E, then  $X_w$ 
is singular if and only if $P_{\idperm,w}(1) > 1$ (Deodhar
\cite{De} proved this for type A, while Peterson (unpublished)
proved it for all simply laced groups.  See \cite{CK}).

 Our main result (Theorem \ref{main theorem}) is a lower bound 
for $P_{x,w}(1)$ in terms of Kazhdan-Lusztig polynomials 
of patterns appearing in $x$ and other elements of $W$ determined by
$x$ and $w$.  Here a pattern of an element of $W$ is its image
under a function $\phi\colon W \to W'$, which we define for
any finite Coxeter group and any (not necessarily standard) 
parabolic subgroup $W' \subset W$.  It agrees with the
standard definition of patterns in type A, 
but is more general than the one using signed permutations used in
\cite{Bi} for types B and D. 

One consequence of our result is the following:

\begin{thm} \label{pattern monotonicity} For any parabolic $W' \subset W$, we
have
$P_{\idperm,w}(1) \ge P_{\idperm,\phi(w)}(1)$.
\end{thm} 
\noindent In particular, this gives another proof that $X_{\tilde{w}}$
singular implies $X_w$ singular in type A.  See also the remark after
Theorem \ref{patterns and geometry}.




The definition of the pattern map $\phi$ is combinatorial, but
it is motivated by the geometry of the action of the torus $T$ on
$\flag$, and the proof of Theorem \ref{main theorem} is entirely 
geometrical.  For $W'\subset W$ parabolic, there is a cocharacter
 $\rho\colon \C^*\to T$ whose fixed point set in $\flag$ 
is a disjoint union of copies of  
the flag variety $\cF'$ of a group $G'$ with Weyl group $W'$.
The action of $\rho$ gives rise to a 
``hyperbolic localization'' functor which takes
sheaves on $\flag$ to sheaves on $\cF'$.  Theorem \ref{main theorem} 
then follows from a ``decomposition theorem'' for this functor, 
proved in \cite{Br},
together with the fact 
that hyperbolic localization preserves local Euler characteristics.

If the action is totally attracting or repelling 
near a fixed point, 
hyperbolic localization is just ordinary restriction or its
Verdier dual.  This gives stronger coefficient-by-coefficient inequalities  
in some special cases (see Theorem \ref{parabolic main}).  The 
attracting/repelling case of \cite{Br} 
has been known for some time; it was used 
in \cite{BrM} to prove a conjecture of Kalai 
on toric $g$-numbers of rational convex polytopes.

Matthew Dyer has recently given us a preprint \cite{Dy} containing an
inequality equivalent to Theorem \ref{main theorem}, which he proves using
his theory of abstract highest weight categories.  It seems likely that
his approach is dual to ours under some version of Koszul duality
\cite{BGS}.

This work was originally motivated by the following question
asked by Francesco Brenti: How can we describe the
Weyl group elements $w$ such that $P_{\id, w}(1)=2$?  In type A, we
can show that if $P_{\id, w}(1)=2$ then the singular locus of the
Schubert variety $X_w$ has only one irreducible component and $w$ must
avoid the patterns: 
$$
\begin{array}{ccc}
(5 2 6 4 1 3)& (5 4 6 2 1 3)& (4 6 3 1 5 2)\\ (4 6 5 1 3 2) &
(6 3 2 5 4 1) & (6 5 3 4 2 1)
\end{array}
$$
We conjecture the converse holds as well.  

We outline the sections of this paper.  In \S\ref{classical.patterns}, we
discuss pattern avoidance on permutations and some applications from
the literature.  In \S\ref{patterns} we describe the pattern map for
arbitrary finite Coxeter groups.
\S\ref{classical and Coxeter} explains why the two notions agree for
permutations.  The main result of \S\ref{patterns} is proved in
\S\ref{coset proof}.  In \S\ref{main result} we state our main
theorem.  In \S\ref{applications} we highlight two particularly
interesting special cases, including Theorem \ref{pattern
monotonicity}.  Our geometric arguments are in \S\ref{geometry}.

\section{Pattern avoidance}
\subsection{Classical pattern avoidance}\label{classical.patterns}

We can write an element $w$ of the permutation group 
${\mfS}_{n}$ on $n$ letters in one-line
notation as $w=w_{1}w_{2}\dotsb w_{n}$, i.e.\ $w$ maps $i$ to $w_{i}$.
We say a permutation $w$ \textit{contains} a pattern $v \in \mfS_{k}$ if
there exists a subsequence $w_{i_{1}}w_{i_{2}}\dotsb w_{i_{k}}$, with
the same relative order as $v=v_{1}\dotsb v_{k}$.  If no such
subsequence exists we say $w$ \textit{avoids} the pattern $v$.

More formally, let $a_{1}\dotsb a_{k}$ be any list of distinct
positive integers.  Define the \textit{flattening function}
$\fl(a_{1}\dotsb a_{k})$ to be the unique permutation $v\in \mfS_{k}$
such that $v_{i}>v_{j} \ \iff \ a_{i}>a_{j}$.  Then it is equivalent
to say that $w$ \textit{avoids} $v$ if no $\fl(w_{i_{1}}w_{i_{2}}\dotsb
w_{i_{k}})=v$.  For example, $w=4 5 3 6 1 7 2$ contains the pattern
$3412$, since $\fl(w_1 w_{4} w_{5} w_{7}) = \fl(4612) = 3412$, but it
avoids $4321$.

Several properties of permutations have been characterized by pattern
avoidance and containment.  For example, as mentioned in the introduction,
for the Schubert variety $X_w$ we have
Schubert variety $X_{w}$ is nonsingular if and only if $P_{\idperm,w}=1$ 
if and only if 
 $w$ avoids $3412$ and $4231$ \cite{LS,Car,De,KL2}.
The element $C_{w}'$ of the Kazhdan-Lusztig basis of the Hecke
algebra of $W$ equals the product $C'_{s_{a_{1}}}
C'_{s_{a_{2}}}\dotsb C'_{s_{a_{p}}}$ for any reduced expression
$w = s_{a_{1}}s_{a_{2}}\dotsb s_{a_{p}}$ if and only if $w$ is
\textit{321-hexagon-avoiding} \cite{BiW}.  Here 321-hexagon-avoiding
means $w$ avoids the five patterns $321,$ $56781234,$ $46781235,$
$56718234,$ $46718235$.

The notion of pattern avoidance easily generalizes to the Weyl groups
of types B,C,D since elements can be represented in one-line notation
as permutations with $\pm$ signs on the entries.  Once
again, the properties $P_{\idperm,w}=1$ and $C_{w}' = C'_{s_{a_{1}}}
C'_{s_{a_{2}}}\dotsb C'_{s_{a_{p}}}$ can be characterized by pattern
avoidance \cite{Bi,BiW}, though the list of patterns can be rather long.
More examples of pattern avoidance appear in
\cite{LasSc,Stem,BiP,BiW2,Manivel,KLR,Co,Co2}.

\subsection{Patterns in Coxeter groups}\label{patterns}
In this section, we generalize the flattening function for permutations
to an arbitrary finite Coxeter group $W$.

Let $S$ be the set of simple
reflections generating $W$.  The set $R$ of all reflections
is $R = \bigcup_{w\in W} wSw^{-1}$.  Given $w\in W$, its
length $l(w)$ is the length of the shortest expression for $w$ 
in terms of elements of $S$.  
The Bruhat-Chevalley order is the partial order $\le$ on $W$ 
generated by the relation
\[x < y \text{  if  } l(x) < l(y)\;\text{and}\; xy^{-1} \in R.\]

Each subset $I \subset S$ generates a subgroup $W_I$; a subgroup 
$W'\subset W$ which is conjugate to $W_I$ for some $I$ is called
a {\em parabolic} subgroup.  The $W_I$'s themselves are known as
{\em standard} parabolic subgroups.

A parabolic subgroup $W' = xW_Ix^{-1}$ of $W$ is again a Coxeter group, 
with simple reflections $S' = xIx^{-1}$ and reflections $R' = R \cap
W'$.  Note that $S' \not\subset S$ unless $W'$ is standard.

We denote the length function and the Bruhat-Chevalley order for $(W',S')$
by $l'$ and $\le'$, respectively.  If $W' = W_I$
then 
\[l' = l|_{W'}\;\text{and}\; \mathord{\le'} = \mathord{\le}|_{W'\times W'},\]
but in general we only have $l'(w) \le l(w)$ and $x \le' y \implies x \le y$.
For instance, if $W' \subset \mfS_4$ is generated by the reflections
$r_{23} = 1324$ and $r_{14} = 4231$, then $r_{23} \le r_{14}$ although they are
not comparable for $\le'$.

The following theorem/definition generalizes the flattening function for 
permutations.

\begin{thm} \label{coset theorem} Let $W' \subset W$ be a parabolic subgroup.
There is a unique function $\phi\colon W\to W'$, the 
{\em pattern map} for $W'$, satisfying:
\begin{enumerate} \item[(a)] 
$\phi$ is $W'$-equivariant: $\phi(wx) = w\phi(x)$ for all $w\in W'$, $x\in W$,
\item[(b)] If $\phi(x) \le' \phi(wx)$ for some $w\in W'$, 
then $x \le wx$.
\end{enumerate}
In particular, $\phi$ restricts to the identity map on $W'$.
\end{thm}
If $W' = W_I$ is a standard parabolic, then (b) can be strengthened to
``if and only if''.  In this case the result is well-known.

To show uniqueness, note that (a) implies that 
$\phi$ is determined by the set $\phi^{-1}(1)$, 
and (b) implies that $\phi^{-1}(1) \cap W'x$ is 
the unique minimal element in $W'x$.  Existence is more subtle; it is
not immediately obvious that the function so defined satisfies (b).
We give a construction of a function $\phi$ that satisfies (a) and 
(b) in Section \ref{coset proof}.  

\subsection{Relation with classical patterns} \label{classical and Coxeter}

Take integers $1 \le a_1 < \dots < a_k \le n$, and let 
$\Sigma = \{a_1, a_2,\ldots, a_k\}$.
Define a generalized
flattening function $\fl_\Sigma\colon \mfS_n \to \mfS_k$
by $\fl_\Sigma(w)=\fl(w_{i_1}w_{i_2}\ldots w_{i_k})$,
 where $w_{i_j} \in \Sigma$ for all $1 \leq j \leq k$ and 
$1 \leq i_1 < i_2 < \ldots < i_k \leq n$.  

Let $W'\subset\mfS_n$ be the subgroup
generated by the transpositions $r_{a_i, a_j}$ for all $i<j$.  It is
parabolic; conjugating by any permutation $z$ with $z_i = a_i$ for
$1 \le i \le k$ gives an isomorphism $\iota\colon \mfS_k \to W'$,
where $\mfS_k \subset \mfS_n$ consists of 
permutations fixing the elements $k+1,\dots, n$.

The function $\iota\circ \fl_\Sigma$ 
satisfies the properties of Theorem~\ref{coset theorem}, 
and so $\iota\circ \fl_\Sigma(w) = \phi(w)$.  Property (a)
follows since left multiplication by a
permutation $w \in W'$ acts only on the values in the set $\{a_1,
a_2,\ldots, a_k\}$.  To prove (b), note that if $v_i = w_i$
for two permutations $v,w \in \mfS_n$ then $v \leq w$ if and only
$\fl(\hat v)\leq \fl(\hat w)$ where $\hat v,\hat w$ 
are the sequences obtained by
removing the $i$th entry from each.
This implies that 
$\iota\circ\fl_\Sigma(x) \leq' \iota\circ \fl_\Sigma(wx)$
if and only if $x \leq wx$.  

For example, take $\Sigma = \{1,4,6,7\}$; the associated subgroup 
$W' \subset \mfS_7$ is
generated by $\{r_{14}, r_{46}, r_{67}\}$.  If $x=6
2 1 3 4 7 5$ then $y=1 2 4 3 6 7 5$ is the unique minimal element in
$W'x$ and $x= r_{46}r_{14} y$, so $\phi(x)=r_{46}r_{14}$.
This agrees with the classical flattening using the isomorphism
$W' \cong \mfS_4$ given by $r_{14} \mapsto s_1$, $r_{46} \mapsto s_2$,
$r_{67} \mapsto s_3$: in fact,
\[\fl_{\{1,4,6,7\}}(6 2 1 3 4 7 5)=\fl(6147) = 3124=s_{2}s_{1}.\]

To obtain the most general parabolic subgroup of $\mfS_n$, let
$\Sigma_1,\dots,\Sigma_l$ be disjoint subsets of $1\dots n$.  To
each $\Sigma_j$ is associated a parabolic subgroup $W'_j$ as before,
and then 
\[W' = W'_1W'_2\dots W'_l \cong \mfS_{|\Sigma_1|} \times \dots \times
\mfS_{|\Sigma_l|}\]
is a parabolic subgroup.  The corresponding flattening function
is \[w \mapsto (\fl_{\Sigma_1}(w),\dots,\fl_{\Sigma_l}(w)).\] 

In types B and D, the flattening function of \cite{Bi} given in terms of 
signed permutations can also be viewed as an instance of our
pattern map.  The group $W'$ of 
signed permutations which fix every element except possibly the
$\pm a_i$, $1 \le i \le k$ is parabolic.  
Multiplication on the left by $w \in W'$ acts only on the
values in the set $\{\pm a_1, \pm a_2,\ldots, \pm a_k\}$ 
and if $v_i = w_i$ for two signed permutations $v,w$ then 
$v \leq w$ if and only $\fl(\hat v)\leq
\fl(\hat w)$ where $\hat v,\hat w$ are the sequences obtained by removing the
$i$th entry from each.  It follows that $v \mapsto \fl(\hat v)$ 
satisfies the conditions of Theorem \ref{coset theorem}     

There are other types of parabolic subgroups in types B and D which 
give rise to other pattern maps.
For instance, the group $W'$ of all unsigned permutations is 
a parabolic subgroup of either B${}_n$ or D${}_n$.  
In this case the pattern map ``flattens'' the signed permutation to an 
unsigned one (e.g. $-4,2,1,-3 \mapsto 1432$). 
Other cases of pattern maps for classical groups are more
difficult to describe combinatorially.

The first author and Postnikov \cite{BiP} have used these more general pattern
maps to reduce significantly the number of patterns needed to recognize
smoothness and rational smoothness of Schubert varieties.  
They reduce the list even further 
by generalizing pattern maps to the case of ``root system
embeddings'' which do not necessarily
preserve the inner products of the roots; for instance, there is a 
root system embedding of A${}_3$ into B${}_3$.  
We do not know of a geometric interpretation of these
more general pattern maps.

\subsection{Spanning subgroups and the reflection representation} 
\label{coset proof}
To prove Theorem \ref{coset theorem} we
use the action of $W$ on its root system.  See \cite[Section 1]{H} 
for proofs of the following facts.

We have the following data: a representation of $W$ on a 
finite-dimensional real vector space $V$, a $W$-invariant subset 
$\Phi \subset V$ (the roots), a subset $\Pi \subset \Phi$ (the positive
roots), and a bijection $r \mapsto \alpha_r$ between  $R$ and  $\Pi$.

These data satisfy the following properties:  $\Phi$ is the disjoint
union of $\Pi$ and $-\Pi$.
The vectors $\{\alpha_s\}_{s\in S}$ form a basis for $V$; a root
$\alpha \in \Phi$ is positive if and only if it  
can be expressed in this basis with nonnegative 
coefficients.  For any $r\in R$ and $w \in W$, we have 
\begin{equation}\label{order and roots}
rw > w \iff \alpha_r \in w\Pi.
\end{equation}

Given a linear function $H\colon V\to \R$, define
\[\Pi_H = \{\alpha \in \Phi\mid H(\alpha)>0\}.\]  
Call $H$ {\em generic\/} if $\Phi \cap\ker H  = \emptyset$.
If we take $H_1(\alpha_s) = 1$ for all $s\in S$, then
$H_1$ is generic and $\Pi = \Pi_{H_1}$.
If we put $H_w = {H_1 \circ w^{-1}}$, then $\Pi_{H_w} = w\Pi$.
Conversely, if $H$ is generic, then $\Pi_H = w\Pi$ for a unique $w\in W$. 

\begin{prop} \label{geometric spanning} Let $W'\subset W$ be a subgroup 
generated by reflections.  Then $W'$ is parabolic if and only if  
there is a subspace $V'\subset V$ so that 
$W'$ is generated by $R' = \{r\in R\mid \alpha_r \in V'\}$.
If so, then $V'$ is $W'$-stable, and putting 
$\Phi' = \Phi\cap V'$, $\Pi' = \Pi \cap V'$,
and $\alpha'_r = \alpha_r$ for $r \in R'$ gives the
reflection representation of $W'$.
\end{prop}

\begin{proof} See \cite[\S1.12]{H}.
\end{proof} 

\begin{rmk}
In type A, all subgroups generated by reflections are parabolic.  
In other types this is no longer the case -- for instance, the subgroup
$W' \cong (\Z_2)^n$ of B${}_n$ generated by reflections in the
roots $\{\pm e_j\}$ is not parabolic for any $n \ge 2$, since
these roots span $V$.
\end{rmk}

We now prove the existence of the function $\phi$ from 
Theorem \ref{coset theorem}.  Let $V' \subset V$ be as 
in Proposition \ref{geometric spanning}. Given $w \in W$,  
we have $w\Pi = \Pi_{H_w}$, and so $w\Pi\cap V' = \Pi'_{H'}$,
where $\Pi' = \Pi\cap V$ and $H' = H_w|_{V'}$.  
It follows that there is a unique $\phi(w) \in W'$ so that
\begin{equation*}\phi(w)\Pi' = w\Pi\cap V'.\end{equation*}
We show that the function $\phi$ defined this way satisfies (a)
and (b) from Theorem \ref{coset theorem}.  
Any $w\in W'$ fixes $V'$, so if $x \in W$ then 
\[\phi(wx) \Pi' = 
(wx\Pi) \cap V' = w(x\Pi \cap V') = w\phi(x)\Pi',\]  
giving (a).

To prove (b), it will be enough to show that 
$\phi(x) \le' \phi(rx)$ implies $x \le rx$ for any $x\in W$, $r\in R'$,
since these relations generate the Bruhat-Chevalley orders on $W$ and $W'$.
We have
\begin{eqnarray*} \phi (x)<'\phi (rx)=r\phi(x) &\iff& \alpha_{r}
\in \phi (x)\Pi'=x\Pi \cap V'\\ & \implies & \alpha_{r} \in x\Pi \\
&\implies& x <rx.
\end{eqnarray*}

\section{The main result}\label{main result}

Suppose now that $W$ is the Weyl group of a semisimple complex algebraic
group $G$.  Let $W' \subset W$ be parabolic, and let $\phi\colon
W\to W'$ be the pattern map of Theorem~\ref{coset
theorem}.  For any $x\in W$, define a partial order on $W'x$ by ``pulling
back'' the Bruhat order from $W'$: if $w,w' \in W'$, say
$wx \le_x w'x$ if and only if $\phi(wx) \le' \phi(w'x)$.
By Theorem \ref{coset theorem}, this is weaker than the Bruhat order on 
$W'x$.

Our main result is the following.

\begin{thm} \label{main theorem} If $x,w \in W$, then
\begin{equation*} 
P_{x,w}(1) \ge \sum_{y\in M(x,w;W')} P_{y,w}(1)P'_{\phi(x), \phi(y)}(1),
\end{equation*}
where $M(x,w;W')$ is the set of maximal elements with respect to 
$\le_x$ in $[1, w] \cap W'x$,
and $P'$ denotes the Kazhdan-Lusztig polynomial 
for the Coxeter system $(W',S')$.
\end{thm}

Conjecturally this should hold for any finite Coxeter group $W$.  There is
a stronger formulation when $W'$ is a standard parabolic subgroup of $W$; 
see the next section.

\begin{ex} Take $W = \mfS_4$, $w = 4231$, $x = 2143$.
Let $W' \cong \mfS_2 \times \mfS_2$ be the group 
generated by reflections $r_{13} = 3214$, $r_{24} = 1432$.
Then $W'x = \{2143, 4123, 2341, 4321\}$.  All but $4321$ 
are in the interval $[1,w]$,
so the maximal elements of $[1,w] \cap W'x$ are $4123 = r_{24}x$ and 
$2341 = r_{13}x$.
Theorem \ref{main theorem} gives
\[P_{2143,4231}(1) \ge P_{4123,4231}(1)P'_{\idperm,r_{24}}(1) + P_{2341,4231}(1)P'_{\idperm,r_{13}}(1) \]
\[ = 1\cdot 1 + 1 \cdot 1 = 2,\]
which holds since $P_{2143,4231}(q) = 1 + q$.

Note that this shows $X_{4231}$ is singular, even though all the Schubert varieties
corresponding to terms on the right hand side are smooth.
\end{ex}


\begin{ex}
One can calculate $P_{1234567,6734512}(1)=44$ in type $A$.  This is
the maximum value of $P_{x,w}(1)$ for any $x,w \in \mfS_{7}$.  
Let $W' \subset \mfS_9$ be the subgroup
generated by the reflections $\{r_{13},r_{34},r_{45},r_{57},r_{78},r_{89}\}$;
it is a parabolic subgroup isomorphic to $\mfS_{7}$.  If 
$w=869457213$ and $x=163457289$,
then $W'x=W'w$ so $M(x,w;W')=\{w \}$, giving $\phi(x)=1234567$ and
$\phi(w)=6734512$.  Hence 
\[P_{x,w}(1)\geq 
P'_{1234567,6734512}(1)P_{w,w}(1)=44.\]
\end{ex}

\subsection{Special cases/applications}\label{applications}
The complicated interaction of the multiplicative structure of $W$ and the
Bruhat-Chevalley order makes computing the set $M(x,w;W')$  difficult.
We mention two cases in which the answer is nice:

\vspace{.1in}
\noindent \textbf{(a)} 
If $w$ and $x$ lie in the same $W'$-coset 
then $M(x,w;W') = \{w\}$.  In this case Theorem \ref{main theorem} says 
\[P_{x,w}(1) \ge P'_{\phi(x), \phi(w)}(1).\]

This allows us to prove Theorem \ref{pattern monotonicity} from the 
introduction': given $w\in W$, let $x \in W'w$ satisfy $\phi(x)= 1$.
Then
\[P_{\idperm,w}(1) \ge P_{x,w}(1) \ge P'_{\idperm,\phi(w)}(1).\]
The first inequality comes from the monotonicity of Kazhdan-Lusztig polynomials
\cite{I},\cite[Corollary 3.7]{BrM2}.

\vspace{.1in}
\noindent \textbf{(b)} 
If either $W'$ or $x^{-1}W'x$ is a standard parabolic subgroup of $W$, then
$M(x,w;W')$ has only one element.  The case where 
$x = 1$ was studied by Billey, Fan, and Losonczy \cite{BiFL}.

In this case the inequality will hold coefficient by coefficient 
rather than just at $q = 1$: 
\begin{thm} \label{parabolic main} If $W'$ or $x^{-1}W'x$ is a standard
parabolic subgroup, then
\[[q^k]P_{x,w} \ge \sum_{i+j = k} [q^i]P_{y,w}[q^j]P'_{\phi(x),\phi(y)},\]
where $M(x,w;W') = \{y\}$.  
Here the notation $[q^k]P$ means the coefficient of $q^k$ in the 
polynomial $P$.
\end{thm}

If both (a) and (b) hold, then Theorem \ref{parabolic main}
 is implied by a well-known equality (see \cite[Lemma 2.6]{P}):
\begin{thm}\label{t:parabolic}
If $W'$ or $x^{-1}W'x$ is a standard parabolic subgroup of $W$ and 
$w \in W'x$, then 
\[P_{x,w}(q) = P'_{\phi(x),\phi(w)}(q).\]
\end{thm}

Theorem~\ref{t:parabolic} can be thought of as a generalization of a
theorem due to Brenti and Simion:
\begin{thm} \cite{BreS} \label{t:BS}
Let $u,v \in \mfS_{n}$.  For any $1 \leq i\leq n$ such
that $\{1,2,\dots,i \}$ appear in the same set of positions (though
not necessarily in the same order) in both $u$ and $v$, then 
\[
P_{u,v}(q)= P_{u[1,i],v[1,i]}(q) \cdot P_{\fl(u[i+1,n]),\fl(v[i+1,n]) }(q),\]
where $u[j,k]$ is obtained from $u$ by only keeping
the numbers $j,j+1,\dots, k$ in the order they appear in $u$.  
\end{thm}

We demonstrate the relationship between the two theorems on an
example.  Let $I_1 = \{s_{1},s_{2},s_{3}\}$, $I_2 = \{s_{5},s_{6},s_{7}\}$,
$I = I_1 \cup I_2$.  Let $W' = W_I \cong W_{I_1} \times W_{I_2}$.
Any pair $x,w$ in the same coset of $W'\backslash W$
satisfies the conditions of
Theorem~\ref{t:BS} and Theorem~\ref{t:parabolic}.  Take $x=25174683$ and
$w=48273561$.  Then Theorem~\ref{t:parabolic} gives
\begin{align*}\label{e:BS.1}
P_{25174683,48273561}(q) &= P'_{\phi(25174683),\phi(48273561)}(q) \\
&= P'_{21435768,42318756}(q) =  P_{2143,4231}(q) P_{1324,4312}(q) 
\end{align*}
agreeing with Theorem~\ref{t:BS}.  The last equality results because
we have $P_{x_1\times x_2, w_1 \times w_2}(q) = P_{x_1,w_1}(q)P_{x_2,w_2}(q)$
for any $x_1 \times x_2$, $w_1 \times w_2$ in the reducible
Coxeter group $W_{I_1}\times W_{I_2}$.

\section{Geometry of flag varieties} \label{geometry}

Let $G$ be a connected 
semisimple linear algebraic group over $\C$.  It acts transitively 
on the flag variety $\flag$ of Borel subgroups of $G$ by conjugation:  
$g\cdot B = gBg^{-1}$.  For any $g \in G$, the point $B\in \cF$
is fixed by $g$ if and only if $g \in B$.

Fix a Borel subgroup and a maximal torus $T \subset B \subset G$.
The Weyl group $W = N_G(T)/T$ is a finite Coxeter group.
The point $g\cdot B \in \cF$ is fixed by $T$ if and only if 
$g \in N_G(T)B$, and so 
$g\mapsto g\cdot B$ induces a bijection between $W$
and $\flag^T$.  We abuse notation and 
refer to $w \in W$ and the corresponding point 
of $\flag$ by the same symbol.

Every $B$-orbit on $\cF$ contains a unique $T$-fixed point;
for $w \in W$, the Bruhat cell $C_w$ is the $B$-orbit $B\cdot w$.
The Schubert variety $X_w$ is the
closure of $C_w$; we have $X_w = \bigcup_{x\le w} C_x$ and
so $X_x \subset X_w \iff x \le w$.

\subsection{Torus actions}\label{T actions}
Let $\rho\colon \C^* \to T$ be a cocharacter of $T$, 
and let $G'$ be the centralizer of $T_0 = \rho(\C^*)$.

\begin{thm} \cite[Theorem 6.4.7]{Sp}\label{restricting Borels} $G'$ 
is connected and reductive; $T$ is a maximal torus in $G'$.  
If $T_0$ fixes a point $B_0 \in \cF$, so that $T_0 \subset B_0$, 
then $B_0 \cap G'$ is a Borel subgroup of $G'$. 
\end{thm}
 Let $\flag'\cong G'/B'$ 
be the flag variety of $G'$, and put $\flag^\rho = \flag^{T_0}$. 
Using Theorem \ref{restricting Borels}, we can define a $G'$-equivariant
algebraic map $\psi\colon \flag^{\rho}\to \flag'$ by 
$\psi(B_0) = (B_0) \cap G'$.

Fix a maximal torus and Borel subgroup of $G'$ by setting 
$B' = B\cap G'$, $T' = T$. 
The Weyl group of $G'$ is $W' = N_{G'}(T')/T' = W \cap (G'/B')$.
The Schubert varieties of $\cF'$ defined by the action of $B'$ 
are indexed by elements of $W'$; denote them by $X'_w$, $w\in W'$.

\begin{prop} \label{tori exist} $W'$ is a parabolic subgroup of $W$, and all
parabolic subgroups arise in this way for some choice of $\rho$.
\end{prop}

This is well-known; the groups $G'$ which arise this way are
Levi subgroups of parabolic subgroups of $G$.  The second half of
the statement (which is the only part we need) can be deduced from
\cite[6.4.3 and 8.4.1]{Sp}, for instance.



Now we can connect the pattern map
$\phi$ defined by Theorem \ref{coset theorem} to geometry. 

\begin{thm} \label{patterns and geometry}
The map $\psi$ 
restricts to an isomorphism on each connected component of $\flag^{\rho}$.
The restriction $\psi|_{\cF^T} \colon \cF^T \to (\cF')^T$ 
is the pattern map $\phi$, using the identifications $\cF^T = W$,
$(\cF')^T = W'$.  In particular, the components of $\cF^\rho$ are in
bijection with $W'\backslash W$.
\end{thm}

\begin{proof} To show the first assertion, it is enough to show that 
$\psi$ is a finite map, since it is $G'$-equivariant and its 
image $\cF'$ is maximal among the compact homogeneous spaces for
$G'$.  But $\psi(g\cdot B) \in (\cF')^T \implies T\subset g\cdot B
\implies g\cdot B \in \cF^T$, a finite set.


Certainly $\psi$ takes $T$-fixed points to $T$-fixed points, so 
it induces a function $W \to W'$ by restriction.
We need to show that it
satisfies the properties of Theorem \ref{coset theorem}.  The
$W'$-equivariance (a) follows immediately from the $G'$-equivariance
of $\psi$.

To see property (b), take $x \in W$ and $w \in W'$, and suppose that
$\psi(x) \le' \psi(wx)$.  This implies that 
$\psi(x) \in \ol{B'\cdot \psi(wx)}$,
and since $x$ and $wx$ lie in the same component of $\flag^{\rho}$, 
we must have
$x \in \ol{B'\cdot wx} \subset \ol{B \cdot wx}$. Thus  
$x \le wx$.   
\end{proof} 

\begin{rmk}\label{r:t-fixed-pts}
Given $w\in W$, let $Y \cong \flag'$ be the component 
of $\flag^{\rho}$ which contains $w$.  Then one can show that 
$X_w \cap Y \cong X'_{\phi(w)}$. Therefore, $X'_{\phi(w)}$ singular 
implies that $X_w$ is singular, using the
result of Fogarty and Norman \cite{FN}: a linearly algebraic group $G$ is 
linearly reductive (this class includes all tori) if and only if for all 
smooth algebraic $G$-schemes $X$ the fixed point scheme $X^{G}$ is smooth.
\end{rmk}

\subsection{Hyperbolic localization}\label{T localization}
Let $X$ be a normal complex variety with an action of
$\C^*$. Let  $X^\circ = X^{\C^*}$, and let $X^\circ_1 \dots X^\circ_r$ be the connected 
components of $X^\circ$. 
For $1 \le k \le r$, define a variety
\[X^+_k = \{x \in X \mid \lim_{t\to 0} t\cdot x \in X^\circ_k\}, 
\]
and let $X^+$ be the disjoint (disconnected) union of all the $X^+_k$.
The inclusions $X^\circ_k \subset X^+_k \subset X$ induce maps
\[X^\circ \ra^f X^+ \ra^g X.\]

Let $D^b(X)$ denote the constructible derived category of 
$\Q$-sheaves on $X$.  
\begin{defn} Given $S\in D^b(X)$, define its
{\em hyperbolic localization\/} 
\[S^{!*} = f^!g^*S \in D^b(X^\circ).\]
\end{defn}

Hyperbolic localization is better adapted to 
$\C^*$-equivariant geometry than ordinary restriction.  It was
first studied by Kirwan \cite{Ki}, who showed that 
if $S$ is the intersection cohomology sheaf of a projective variety
with a linear $\C^*$-action,
then $S$ and $S^{!*}$ have isomorphic hypercohomology groups.  

We will need two properties of hyperbolic localization
from \cite{Br}.
For any $S \in D^b(X)$ and $p \in X$, we let $\chi_p(X)$ 
denote the Euler characteristic of the stalk cohomology at $p$.

\begin{prop}\cite[Proposition 3]{Br} 
\label{localization chi} If $p \in X^\circ$, then 
\[\chi_p(S) = \chi_p(S^{!*}).\]
\end{prop}

Second, hyperbolic localization satisfies a  
decomposition theorem \cite[Theorem 2]{Br}. When applied to 
$X = \flag$ and the action given by $\rho$, this gives the following.

\begin{thm} \label{localization purity}
Let $L_w$, and $L'_v$ be the intersection
cohomology sheaves of the Schubert varieties $X_w$ and $X'_v$, 
respectively.
For any $w\in W$ and $1 \le k \le r$, there is an isomorphism
\begin{equation*}
\psi_*((L_w)^{!*}|_{\cF^\circ_k}) \cong \bigoplus_{j=1}^{m} L'_{v_j}[d_j],
\end{equation*}
for some $v_j \in W'$ (not necessarily distinct) and $d_j \in 2\Z$.
\end{thm}
Here we use the fact that hyperbolic localization preserves 
$B'$-equivariance.  The fact that $d_j \in 2\Z$ follows from the purity 
of the stalks of simple mixed Hodge modules of Schubert
varieties.

\subsection{Proof of Theorem \ref{main theorem}}
The description of Kazhdan-Lusztig polynomials as the 
local intersection cohomology Poincar\'e polynomials
of Schubert varieties \cite{KL2} implies that for any $u, v \in W$, we have
\[P_{u,v}(1) = \chi_u(L_v) = \sum_i \dim_\Q \HH^{2i}((L_v)_u).\]

Now, given $x,w\in W$, let $\cF^\circ_k$ be the component
of $\cF^{\rho}$ which contains $x$, and thus all of $W'x$.  
For every $y\in W'$, let $a_y$
be the number of $j$ for which $v_j = y$ in Theorem 
\ref{localization purity}.  

For any $z \in W'x$ we have, 
using Theorem \ref{localization purity} and Proposition
\ref{localization chi},
\begin{align} \nonumber
P_{z,w}(1)&= \chi_z(L_w)= \chi_{\phi(z)}
\left(\psi_*((L_w)^{!*}|_{\cF^\circ_k})  \right)\\
\label{loc sum} &=\sum_{j=1}^{m} \chi_{\phi(z)}\left(L'_{v_j}[d_j] \right)\\
\nonumber &=\sum_{y\in W'z} a_y P'_{\phi(z),\phi(y)}(1)
\end{align}
(note that the shift $[d_j]$ does not change the
Euler characteristic, since $d_j \in 2\Z$).

If $z \notin [1,w]$ then equation 
(\ref{loc sum}) implies $a_z = 0$,
since $P_{z,w} = 0$, $P'_{z,z} = 1$, and all the terms in the sum
are nonnegative.
Using (\ref{loc sum}) again shows that if $y \in M(x;w;W')$, i.e.\
$y$ is maximal in $[1,w]\cap W'x$, then 
$a_y = P_{y,w}(1)$.  Finally, 
evaluating \eqref{loc sum} at $x$ and keeping only
the terms with $y \in M(x,w;W')$ proves Theorem \ref{main theorem}.
\qed


\subsection{Proof of Theorem~\ref{parabolic main}} 
Suppose first that $x^{-1}W'x = W_I$ is a standard parabolic subgroup.
Take $\mu$ to be any dominant integral cocharacter
which annihilates a root $\alpha_r$ if and only if $r \in W'$,
and let $\rho = Ad(x)\mu$.  Then the action 
of $\rho$ is {\em completely
repelling\/} near the component $\flag^\circ_k$ of $\flag^{\rho}$ 
which contains $W'x = xW_I$, meaning that 
$\flag_k^+ = \flag^\circ_k$, in the notation of \S\ref{T localization}.
 
This implies that hyperbolic localization to 
$\cF^\circ_k$ is just ordinary restriction:
setting $h\colon \cF^\circ_k \to \cF^{\rho}$ for the inclusion, we have
\[(S^{!*})|_{\cF^\circ_k} = h^!f^!g^*S = (fh)^!g^*S = (fh)^*g^*S = 
S|_{\cF^\circ_k},\]  
since both $h$ and $fh$ are open immersions.
The same argument given for Theorem \ref{main theorem} now proves Theorem 
\ref{parabolic main}, using local Poincar\'e polynomials instead of 
local Euler characteristics.  

If instead $W' = W_I$, we can use the anti-involution 
$g \mapsto g^{-1}$  to replace left cosets by right cosets, since
$P_{x^{-1}, w^{-1}} = P_{x,w}$ for all
$x,w \in W$. 

\subsection*{Acknowledgments}

We have benefitted greatly from conversations with Francesco Brenti, 
Victor Guillemin, Victor Ginzburg, Bert Kostant, Sue Tolman, 
David Vogan and Greg Warrington.  We are grateful to Patrick Polo for 
spotting an error in Theorem 5, and to the referees for thoughtful 
comments.


\begin{thebibliography}{99}

\bibitem[BB]{BeBe} A.~Beilinson and J.~Bernstein, {\em 
Localisation de $\mathfrak g$-modules},
C.~R.~Acad.~Sci.~Paris Ser.~I Math. {\bf 292} (1981), no.~1, 15--18. 

\bibitem[BGS]{BGS} A.~Beilinson, V.~Ginzburg, and W.~Soergel 
{\em Koszul duality 
patterns in representation theory}, J.~Amer.~Math.~Soc. {\bf 9} (1996), no.~2, 
473--527.

\bibitem[Bi]{Bi} S.~Billey, {\em Pattern avoidance and rational
smoothness of Schubert varieties} Adv.~Math.~{\bf 139} (1998), no.~1,
141--156.

\bibitem[BiFL]{BiFL} S.~Billey, C.K.~Fan, and J.~Losonczy, {\em The
parabolic map}.  J.~Algebra, {\bf 214} (1999), no.~1, 1--7.

\bibitem[BiP]{BiP} S.~Billey and A.~Postnikov, {\em A root system
description of pattern avoidance with applications to Schubert
varieties},  preprint {\tt math.CO/0205179} (2002), 17 pp. 

\bibitem[BiW]{BiW}  
S.~Billey and G.~Warrington, {\em Kazhdan--Lusztig
Polynomials for 321-hexagon-avoiding permutations}, J.~Alg.~Comb.
\textbf{13} (2000), 111--136.

\bibitem[BiW2]{BiW2} \bysame, {\em Maximal
Singular Loci of Schubert Varieties in $SL_{n}/B$}, 
Trans. Amer. Math. Soc. {\bf 355} (2003),
no.~10, 3915--3945.

\bibitem[Br]{Br} T.~Braden, {\em Hyperbolic localization of intersection
cohomology}, preprint {\tt math.AG/0202251} (2002), 7 pp., submitted.

\bibitem[BrM]{BrM} T.~Braden and R.~MacPherson, {\em Intersection
homology of toric varieties and a conjecture of Kalai},
Comment.~Math.~Helv. {\bf 74} (1999), no.~3, 442--455.

\bibitem[BrM2]{BrM2} \bysame, {\em From moment
graphs to intersection cohomology}, Math.~Ann. {\bf 321} (2001) 
no.~3, 533--551.

\bibitem[BreS]{BreS}
{F.~Brenti and R.~Simion}, 
{\em Explicit formulae for some Kazhdan-Lusztig polynomials}, 
J. Algebraic Combin. {\bf 11} (2000), no. 3, 187--196. 

\bibitem[BryK]{BryK}  J.-L.~Brylinski and M.~Kashiwara, 
{\em Kazhdan-Lusztig conjecture and holonomic systems}, 
Invent.~Math. {\bf 64} (1981), no.~3, 387--410. 

\bibitem[C]{Car} J.~Carrell, {\em The Bruhat graph of a Coxeter
group, a conjecture of Deodhar, and rational smoothness of Schubert
varieties}, Proceedings of Symposia in Pure Math., \textbf{56} (1994),
53--61.

\bibitem[CK]{CK} J.~Carrell and J.~Kuttler, {\em Smooth points of
$T$-stable varieties in $G/B$ and the Peterson map}, to appear in
Invent.~Math. 

\bibitem[Co]{Co} A.~Cortez, 
{\em Singularit\'es g\'en\'eriques des vari\'et\'es de Schubert covexillaires},
Ann.~Inst.~Fourier (Grenoble) {\bf 51} (2001), no.~2, 375-393.

\bibitem[Co2]{Co2} 
A.~Cortez, {\em Singularit\'es g\'en\'eriques et quasi-resolutions 
des vari\'et\'es de Schubert pour le groupe lin\'eaire},
C. R. Acad. Sci. Paris S\'er. I Math. {\bf 333} (2001), no. 6, 561--566. 

\bibitem[De]{De} V.~Deodhar, {\em Local Poincar\'e duality and 
nonsingularity of Schubert varieties},  
Comm.~Algebra {\bf 13} (1985), no.~6, 1379--1388.

\bibitem[Dy]{Dy} M.~Dyer, {\em Rank two detection of singularities of
Schubert varieties}, preprint (2001).

\bibitem[FN]{FN}
{J.~Fogarty and P.~Norman}, {\em A fixed-point characterization of linearly
  reductive groups}, in Contributions to algebra (collection of papers
  dedicated to {Ellis Kolchin}), Academic Press, New York, 1977, pp.~151--155.

\bibitem[H]{H} J.~Humphreys, Reflection groups and Coxeter groups,
Cambridge studies in advanced mathematics {\bf 29}.

\bibitem[I]{I} R.~Irving, {\em The socle filtration of a Verma
module}, Ann.~Sci.~\'Ecole Norm.~Sup. series 4 {\bf 21} (1988), no.~1,
47--65.

\bibitem[KL1]{KL1}  D.~Kazhdan and G.~Lusztig, {\em Representations of
Coxeter groups and Hecke algebras}, Invent.~Math., \textbf{53} (1979),
165--184.

\bibitem[KL2]{KL2} \bysame, {\em Schubert varieties
and Poincar\'e duality}, Proc.~Symp.~Pure.~Math., AMS, \textbf{36}
(1980), 185--203.


\bibitem [Ki]{Ki} F.~Kirwan, {\em Intersection homology and torus
actions}, J.~Amer.~Math.~Soc. {\bf 1} (1988), no.~2, 385--400.

\bibitem[LS]{LS} V.~Lakshmibai and B.~Sandhya, {\em Criterion for
smoothness of Schubert varieties in ${\rm Sl}(n)/B$} Proc.~Indian
Acad.~Sci.~Math.~Sci. {\bf 100} (1990), no.~1, 45--52.

\bibitem[LasSc]{LasSc} A.~Lascoux and M.P.~Sch\"utzenberger, {\em
Polyn\^omes de Kazhdan and Lusztig pour les grassmanniennes.}
(French) [Kazhdan--Lusztig polynomials for Grassmannians],
Ast\'erisque, \textbf{87--88} (1981), 249--266, Young tableaux and
Schur functions in algebra and geometry (Toru\'n, 1980).

\bibitem[KLR]{KLR} C. ~Kassel, A. ~Lascoux  and C.~Reutenauer,
\emph{The singular locus of a {Schubert} variety}, Pr\'epublication de
l'{Institut de Recherche Math\'ematique Avanc\'ee} (2001).  To appear
J.~of Alg.

\bibitem[Ma]{Manivel} L.~Manivel, \emph{Le lieu singulier des
vari\'et\'es de {Schubert}}, Internat. Math. Res. Notices (2001), 
no. 16, 849--871.

\bibitem[P]{P} P.~Polo, 
{\em Construction of arbitrary Kazhdan-Lusztig polynomials in symmetric
groups} Represent.~Theory {\bf 3} (1999), 90-104 (electronic).

\bibitem[R]{R} K. M. Ryan, {\em On Schubert varieties in the flag 
manifold of ${\rm Sl}(n,{\bf C})$}, Math. Ann. {\bf 276} (1987), 
no.~2, 205--224.

\bibitem[Sp]{Sp} T.~Springer, Linear algebraic groups, Second
edition. Birkh\"auser, Boston, MA, 1998.

\bibitem[St]{Stem} J.~Stembridge, {\em On the fully commutative
elements of {Coxeter} groups}, J.~Algebraic Combin. \textbf{5} (1996),
no.~4, 353--385.

\bibitem[W]{W} J. S. Wolper, {\em A combinatorial approach to the 
singularities of Schubert varieties}, Adv. Math. {\bf 76} (1989), 
no.~2, 184--193.

\end{thebibliography}
\end{document}